\documentclass[a4paper,10pt,dvipsnames]{article}
\usepackage[utf8]{inputenc}
\usepackage[english]{babel}
\usepackage{amsthm}
\usepackage{amsfonts}
\usepackage{amssymb}	
\usepackage{amsmath}
\allowdisplaybreaks
\usepackage{chemarr}
\usepackage{slashed}
\usepackage{enumerate}
\usepackage{graphicx}
\usepackage{systeme}
\usepackage{dsfont}
\usepackage{relsize}
\usepackage[margin=0.90in]{geometry}
\usepackage{float}
\usepackage{upgreek}
\usepackage{titlesec}
\usepackage{tikz}
\usepackage[labelformat=simple]{subcaption}

\usepackage{graphicx}
\usepackage{bmpsize}
\usepackage{epstopdf}
\usepackage{bbm}
\usepackage{etoolbox}
\usepackage{xcolor}
\usepackage{titling}

\usepackage{blindtext,graphicx}
\usepackage[absolute]{textpos}
\usepackage[hang,flushmargin]{footmisc}

\usepackage{xfrac}
\usepackage{nicefrac}
\usepackage{empheq}
\usepackage{stmaryrd}
\usepackage{cancel}
\usepackage[linguistics]{forest}
\usepackage{url}
\usepackage[font=small]{caption}
\usepackage{array}
\usepackage{dsfont}
\usepackage{tikz}
\usetikzlibrary{trees}
\usetikzlibrary{calc}
\usepackage{pdflscape}
\usepackage{rotating}

\usepackage{todonotes}
\usepackage{xcolor,hyperref}
\definecolor{darkblue}{rgb}{0,0,0.6}
\hypersetup{
    colorlinks=true,       
    linkcolor=darkblue,          
    citecolor=darkblue,        
    filecolor=darkblue,      
    urlcolor=darkblue           
}

\usepackage[giveninits=true,sortcites=true,date=year,maxbibnames=99,doi=false,isbn=false,url=false,eprint=false]{biblatex}
\renewbibmacro{in:}{}
\DeclareFieldFormat{pages}{#1}
\AtEveryBibitem{%
\clearlist{language}%
}
\usepackage{csquotes}

\usepackage{setspace}
\makeatletter
\newcommand{\thickhline}{%
\noalign {\ifnum 0=`}\fi \hrule height 1pt

\futurelet \reserved@a \@xhline
}
\newcolumntype{"}{@{\hskip\tabcolsep\vrule width 1pt\hskip\tabcolsep}}
\makeatother
\theoremstyle{definition}

\makeatletter
\newcommand{\pushright}[1]{\ifmeasuring@#1\else\omit\hfill$\displaystyle#1$\fi\ignorespaces}
\newcommand{\pushleft}[1]{\ifmeasuring@#1\else\omit$\displaystyle#1$\hfill\fi\ignorespaces}
\makeatother

\titleformat{\paragraph}
{\normalfont\normalsize\bfseries}{\theparagraph}{1em}{}
\titlespacing*{\paragraph}
{0pt}{3.25ex plus 1ex minus .2ex}{1.5ex plus .2ex}

\newcounter{alphasect}
\def\alphainsection{0}

\let\oldsection=\section
\def\section{%
  \ifnum\alphainsection=1%
    \addtocounter{alphasect}{1}
  \fi%
\oldsection}%
\renewcommand\thesection{%
  \ifnum\alphainsection=1%
    \Alph{alphasect}%
  \else%
    \arabic{section}%
  \fi%
}%

\title{The isotropic relaxed micromorphic model in polar coordinates and its application to an elastostatic axisymmetric extension problem}
\author{
Esmaeal Ghavanloo\,\thanks{Corresponding author: Esmaeal Ghavanloo, School of Mechanical Engineering, Shiraz University, Shiraz, 71963-16548, Iran, email: ghavanloo@shirazu.ac.ir}, \quad
 and \quad
Patrizio Neff\,\thanks{Patrizio Neff, Head of Chair for Nonlinear Analysis and Modelling, Fakultät für Mathematik, Universität Duisburg-Essen, Thea-Leymann-Straße 9, 45127 Essen, Germany, email: patrizio.neff@uni-due.de}
}

\thanksmarkseries{arabic}
\date{\today}
\begin{document}
\maketitle
\begin{abstract}
\noindent
In this paper, we consider the isotropic relaxed micromorphic model in polar coordinates and use this representation to solve explicitly an elastostatic axisymmetric extension problem involving a linear system of ordinary differential equations. To obtain an analytical solution, modified Bessel functions are utilized and closed-form solutions for the displacement and microdistortion are obtained. We show how certain limit cases (classical linear elasticity), which are naturally included in the relaxed micromorphic model, can be efficiently achieved. Furthermore, numerical results are calculated and the effects of various parameters are examined. The results can be used to calibrate and check corresponding finite element codes.
\end{abstract}
\textbf{Keywords}: generalized continua, relaxed micromorphic model, polar coordinates; axisymmetric extension problem; consistent coupling
{
%
%
%
%
%
\section{Introduction}
\label{sec:intro}
Experimental observations in both natural and synthetic materials (e.g., metamaterials) demonstrate that their physical behavior varies with size \cite{liu2018experimental,cui2018additive}, making it impossible to accurately predict these behaviors using classical continuum theory. To overcome this limitation, various generalized continuum models have been developed, that extend classical continuum theory for modeling size-dependent responses \cite{bleustein1966effects,ghavanloo2021size}. These generalized continuum models are typically classified into two categories. In the first category, various additional degrees of freedom, such as micro-rotation \cite{lakes2023cosserat}, micro-stretch, and micro-strain \cite{forest2006nonlinear}, are introduced to describe the deformation of the microstructure. The classical micromorphic theory \cite{Mindlin1964,Eringen1999} and the Cosserat theory \cite{papamichos2010continua} are two of the most well-known theories in this category. The second category of generalized continua involves the incorporation of higher-order differential operators in the energy or motion functional \cite{askes2011gradient}. The strain gradient theory is the most well-known theory within that category \cite{alavi2023hierarchy}. 

Among the various generalized continuum models, micromorphic models have been particularly successful in describing the behavior of mechanical metamaterials \cite{alberdi2021exploring}, composites \cite{pouramiri2024estimation}, porous media \cite{ju2024three}, and granular materials \cite{xiu2020micromechanics}. However, the classical micromorphic theory involves a large number of material coefficients (18 in the linear-isotropic case), which restricts its practical use. To address these limitations, various simplified versions of the micromorphic theory have been developed over the past decade \cite{neff2014unifying,romano2016micromorphic,shaat2018reduced,zhang2021non}. The relaxed micromorphic model is a generalized continuum model that effectively describes size effects while significantly reducing the number of material coefficients required compared to the classical theory. Over the past decade, this model has been successfully applied to solve a wide range of static and dynamic problems \cite{gourgiotis2024green,madeo2016complete,owczarek2021note,rizzi2022analytical,rizzi2021analytical2,rizzi2021analytical3,sky2024novel,knees2024global}. The ongoing development and refinement of the relaxed micromorphic model continue to expand its applicability, promising further advancements in the understanding and utilization of complex material systems. 

In the relaxed micromorphic model, the kinematics are defined by the classical displacement $u: \Omega \to \mathbb{R}^3$ and the non-symmetric micro-distortion $P: \Omega \to \mathbb{R}^{3 \times 3}$. The solution is subsequently obtained from the variational two-field problem \cite{sky2022primal}:
{\small
\begin{align}
I(u,P)=
&\int_{\Omega}
\frac{1}{2}
\bigg(
\langle
\mathbb{C}_{\rm e}\,\text{sym} (\text{D}u - P), \text{sym} (\text{D}u - P)
\rangle
+
\langle
\mathbb{C}_{\rm c}\,\text{skew} (\text{D}u - P), \text{skew} (\text{D}u - P)
\rangle
\label{eq:energy_rmm}
\\*
&\phantom{\int_{\Omega} \frac{1}{2} \bigg(}
+
\langle
\mathbb{C}_{\rm micro}\,\text{sym} \, P, \text{sym} \, P
\rangle
+
\mu_{\rm macro}  L_{\rm c}^2
\langle
\mathbb{L}\,\text{Curl} \, P, \text{Curl} \, P
\rangle
\bigg)
\text{d}x
\quad
\longrightarrow
\quad
\text{min}\,(u,P)
\, ,
\notag
\end{align}
}
where $\mathbb{C}_{\rm e},\mathbb{C}_{\rm micro},\mathbb{L}$ are positive-definite fourth-order tensors, $\mathbb{C}_{\rm c}$ is a positive semi-definite fourth order tensor, $L_{\rm c}$ is a characteristic length and $\mu_{\rm macro}$ is the macroscopic shear modulus added for dimensional consistency. In addition, the stiffness tensors $\mathbb{C}_{\rm macro}$, $\mathbb{C}_{\rm e}$ and $\mathbb{C}_{\rm micro}$ are related together by the following relationships \cite{neff2014unifying,barbagallo2017transparent}:
{\small
\begin{gather}
\mathbb{C}_{\rm e}
=
\mathbb{C}_{\rm micro}
\bigg(\mathbb{C}_{\rm micro}-\mathbb{C}_{\rm macro}
\bigg)^{-1}
\mathbb{C}_{\rm macro}
\quad
\Longleftrightarrow
\quad
\mathbb{C}_{\rm macro}
=
\mathbb{C}_{\rm micro}
\bigg(\mathbb{C}_{\rm micro}+\mathbb{C}_{\rm e}
\bigg)^{-1}
\mathbb{C}_{\rm e}
\, ,
\label{eq:Ce_Cmac}
\\
\mathbb{C}_{\rm micro}
=
\mathbb{C}_{\rm e}
\bigg(\mathbb{C}_{\rm e}-\mathbb{C}_{\rm macro}
\bigg)^{-1}
\mathbb{C}_{\rm macro}
\, ,
\notag
\end{gather}
}

The original relaxed micromorphic model, along with its derived models, has been formulated in general tensor forms \cite{demetriou2024reduced}. This formulation allows for theoretical adaptation to specific problems where Cartesian coordinates are suitable. When applying the relaxed micromorphic model in scenarios where curvilinear coordinates (such as polar, cylindrical, or spherical coordinates) are more appropriate, deriving the corresponding formulations for governing equations and boundary conditions is not straightforward. The process is often complex, tedious, and challenging. However, to date, in spite of the vast application of the relaxed micromorphic model in prediction of the mechanical behavior of various natural and man-made materials, its general formulation in orthogonal curvilinear coordinates is absent in the literature. Motivated by this missing development, we derive general formulations for the relaxed micromorphic model in polar coordinates by utilizing the expressions for divergence, curl, and gradient operators of tensors in these coordinates. Then, to demonstrate the practical application of the presented formulation, we analytically solve an elastostatic axisymmetric extension problem. This involves analyzing the deformation of a thin circular plate subjected to uniform radial displacement at its boundary, using modified Bessel's functions. In addition, it is demonstrated that the classical elasticity model can be derived as limit-cases of the relaxed micromorphic solution. Finally, numerical results are calculated to illustrate the effects of various parameters, including material coefficients of the relaxed meromorphic model and the characteristic length.  

%
\section{The isotropic relaxed micromorphic model}
The isotropic relaxed micromorphic model exhibits the same kinematics as the classical micromorphic isotropic model \cite{Mindlin1964,Eringen1999}. In this framework, the displacement $u $ and the non-symmetric microdistortion $P$ are independent fields. The operator form of governing equations of motion in the absence of body forces is given by \cite{madeo2015band}
{\small
\begin{align}
\rho\, \ddot{u}= \text{Div} \, \sigma  \, ,
\qquad\qquad\qquad\qquad
\eta\, \ddot{P}  &=  \sigma - \sigma_{\rm micro}- \text{Curl}\,m\, ,
\label{eq:equi_RM}
\end{align}
where
\begin{align}
\label{constEq}
\sigma & \coloneqq
2\mu_{\rm e}\,\text{sym}  (\text{D} u - P )
+ 2\mu_{\rm c}\,\text{skew} (\text{D} u - P )
+ \lambda_{\rm e} \text{tr} (\text{D} u - P ) \mathbbm{1}
\, ,
\notag
\\*
\sigma_{\rm micro} & \coloneqq
2 \mu_{\rm micro}\,\text{sym}\,P
+ \lambda_{\rm micro} \text{tr} \left( P \right) \mathbbm{1}
\, ,
\\*
\notag
m & \coloneqq
\mu_{\rm macro}\, L_{\rm c}^2 \text{Curl} \, P 
\, ,
\end{align}
}
where $\sigma$ is the non-symmetric elastic force stress tensor, $\sigma_{\rm micro}$ is the symmetric microscopic stress tensor, $m$ is the non-symmetric moment tensor, $\rho$ and $\eta$ are the macro and micro mass densities respectively and all other quantities ($\mu_{\rm e}$, $\lambda_{\rm e}$, $\mu_{\rm c}$, $\mu_{\rm micro}$, $ \lambda_{\rm micro}$, $\mu_{\rm macro}$ and $L_{\rm c}$) are the constitutive parameters of the model. Furthermore, $\mathbbm{1}$ denotes the identity matrix, $\text{Div}$ represents the row-wise divergence operator, $\text{D}$ denotes the gradient operator, $\text{sym}$ and $\text{skew}$ indicate the symmetric and skew-symmetric parts of a tensor, respectively. 

The homogeneous Neumann and the Dirichlet boundary conditions are  
{\small
\begin{align}
&
\text{Neumann:}
&&&&
t
\coloneqq
\sigma \cdot \, n
= 
0
\, ,
&&
\text{and}
&&
\eta
\coloneqq
m \times n = 
0
\, ,
&
\label{eq:BC_RM_Neu}
\\*
&
\text{Dirichlet:}
&&&&
u
=
\overline{u}
\, ,
&&
\text{and}
&&
\overline{Q}=P \cdot \tau 
\, ,
&
\label{eq:BC_RM_Dir}
\end{align}}
in which $\tau$ and $n$ are the outer tangent and normal vector on the boundary. The higher-order Dirichlet boundary conditions in  Eq. (\ref{eq:BC_RM_Dir}) can be particularised to \cite{schroder2022lagrange}
{\small
\begin{align}
P \cdot \tau
=
\overline{Q}
=
\text{D}u \cdot \tau
\, ,
\label{eq:CCBC_RM}
\end{align}
}
called ``consistent coupling boundary conditions'' \cite{d2022consistent}. 

\section{The isotropic relaxed micromorphic model in polar coordinates}
\subsection{Governing equations}
In this section, the formulations of the isotropic two-dimensional relaxed micromorphic model in polar coordinates will be derived. In the polar coordinates, the position of a material point is identified by $r$ and $\theta$, where $r$ is the radial distance from the origin and $\theta$ is the azimuthal angle. The displacement field $u$ and the micro-distortion tensor $P $ can be expressed in the polar coordinates as follows:
{\small
\begin{align}
{u}={{u}_{r}}(r,\theta )\,{{\mathbf{e}}_{\mathbf{r}}}+{{u}_{\theta }}(r,\theta )\,{{\mathbf{e}}_{\mathbf{\theta }}}\, ,
\label{eq:u_polar}
\end{align}

\begin{align}
P={{P}_{rr}}(r,\theta )\,{{\mathbf{e}}_{\mathbf{r}}}\otimes{{\mathbf{e}}_{\mathbf{r}}}+{{P}_{r\theta }}(r,\theta )\,{{\mathbf{e}}_{\mathbf{r}}}\otimes{{\mathbf{e}}_{\mathbf{\theta }}}+{{P}_{\theta r}}(r,\theta )\,{{\mathbf{e}}_{\mathbf{\theta }}}\otimes{{\mathbf{e}}_{\mathbf{r}}}+{{P}_{\theta \theta }}(r,\theta )\,{{\mathbf{e}}_{\mathbf{\theta }}}\otimes{{\mathbf{e}}_{\mathbf{\theta }}}\, ,
\label{eq:P_polar}
\end{align}
}
where ${{\mathbf{e}}_{\mathbf{r}}}$ and ${{\mathbf{e}}_{\mathbf{\theta }}}$ are the unit vectors in $r$ and $\theta$ directions respectively. Furthermore, ${u}_{r}$ and ${u}_{\theta }$ denote the radial and angular displacements, respectively. To present Eqs. (\ref{eq:equi_RM}-\ref{constEq}) in polar coordinates, the gradient of the displacement field, the divergence of the stress tensor, and the curl of the microdistortion tensor must be properly defined in polar coordinates. The gradient of displacement field can be expressed as \cite{sadd2009elasticity}
{\small
\begin{align}
\text{D} {u}=\frac{\partial {{u}_{r}}}{\partial r}\,{{\mathbf{e}}_{\mathbf{r}}}\otimes{{\mathbf{e}}_{\mathbf{r}}}+\frac{1}{r}\left( \frac{\partial {{u}_{r}}}{\partial \theta }-{{u}_{\theta }} \right)\,{{\mathbf{e}}_{\mathbf{r}}}\otimes{{\mathbf{e}}_{\mathbf{\theta }}}+\frac{\partial {{u}_{\theta }}}{\partial r}\,{{\mathbf{e}}_{\mathbf{\theta }}}\otimes{{\mathbf{e}}_{\mathbf{r}}}+\frac{1}{r}\left( \frac{\partial {{u}_{\theta }}}{\partial \theta }+{{u}_{r}} \right)\,{{\mathbf{e}}_{\mathbf{\theta }}}\otimes{{\mathbf{e}}_{\mathbf{\theta }}}\, .
\label{grad_u}
\end{align}
}
In addition, the divergence of the stress tensor $\sigma$ can be express as
{\small
\begin{align}
\text{Div} \, \sigma  =\left( \frac{\partial {{\sigma }_{rr}}}{\partial r}+\frac{1}{r}\frac{\partial {{\sigma }_{r\theta }}}{\partial \theta }+\frac{{{\sigma }_{rr}}-{{\sigma }_{\theta \theta }}}{r} \right){{\mathbf{e}}_{\mathbf{r}}}+\left( \frac{\partial {{\sigma }_{\theta r}}}{\partial r}+\frac{1}{r}\frac{\partial {{\sigma }_{\theta \theta }}}{\partial \theta }+\frac{{{\sigma }_{r\theta }}+{{\sigma }_{\theta r}}}{r} \right){{\mathbf{e}}_{\mathbf{\theta }}}\, ,
\label{Div_sigm}
\end{align}
}
where $\sigma_{rr}$, $\sigma_{r\theta }$, $\sigma_{\theta r}$ and $\sigma_{\theta \theta }$ are the components of the stress tensor.

Now, using Eqs. (\ref{eq:equi_RM}), (\ref{constEq}), (\ref{eq:P_polar})-(\ref{Div_sigm}), we obtain
{\small
\begin{align}\label{Gov_u_r}
\nonumber  \rho \frac{{{\partial }^{2}}{{u}_{r}}}{\partial {{t}^{2}}}=\,&{{\mu }_{\rm e}}\left( 2\left[ \frac{{{\partial }^{2}}{{u}_{r}}}{\partial {{r}^{2}}}+\frac{1}{r}\frac{\partial {{u}_{r}}}{\partial r}-\frac{{{u}_{r}}}{{{r}^{2}}}+\frac{{{P}_{\theta \theta }}-{{P}_{rr}}}{r}-\frac{\partial {{P}_{rr}}}{\partial r} \right]+\frac{1}{{{r}^{2}}}\frac{{{\partial }^{2}}{{u}_{r}}}{\partial {{\theta }^{2}}}-\frac{3}{{{r}^{2}}}\frac{\partial {{u}_{\theta }}}{\partial \theta }+\frac{1}{r}\frac{{{\partial }^{2}}{{u}_{\theta }}}{\partial r\partial \theta }-\frac{1}{r}\frac{\partial {{P}_{r\theta }}}{\partial \theta }-\frac{1}{r}\frac{\partial {{P}_{\theta r}}}{\partial \theta } \right) \\
\nonumber
& +{{\lambda }_{\rm e}}\left( \frac{1}{r}\frac{{{\partial }^{2}}{{u}_{\theta }}}{\partial r\partial \theta }-\frac{\partial {{P}_{\theta \theta }}}{\partial r}+\frac{{{\partial }^{2}}{{u}_{r}}}{\partial {{r}^{2}}}+\frac{1}{r}\frac{\partial {{u}_{r}}}{\partial r}-\frac{{{u}_{r}}}{{{r}^{2}}}-\frac{1}{{{r}^{2}}}\frac{\partial {{u}_{\theta }}}{\partial \theta }-\frac{\partial {{P}_{rr}}}{\partial r} \right) \\ 
& +{{\mu }_{\rm c}}\left( \frac{1}{{{r}^{2}}}\frac{{{\partial }^{2}}{{u}_{r}}}{\partial {{\theta }^{2}}}-\frac{1}{{{r}^{2}}}\frac{\partial {{u}_{\theta }}}{\partial \theta }-\frac{1}{r}\frac{{{\partial }^{2}}{{u}_{\theta }}}{\partial r\partial \theta }-\frac{1}{r}\frac{\partial {{P}_{r\theta }}}{\partial \theta }+\frac{1}{r}\frac{\partial {{P}_{\theta r}}}{\partial \theta } \right),
\end{align}

}
{\small
\begin{align}\label{Gov_u_teta}
 \nonumber  \rho \frac{{{\partial }^{2}}{{u}_{\theta }}}{\partial {{t}^{2}}}=\,&{{\mu }_{\rm e}}\left( \frac{1}{r}\frac{{{\partial }^{2}}{{u}_{r}}}{\partial r\partial \theta }+\frac{1}{{{r}^{2}}}\frac{\partial {{u}_{r}}}{\partial \theta }+\frac{2}{{{r}^{2}}}\frac{\partial {{u}_{r}}}{\partial \theta }+\frac{{{\partial }^{2}}{{u}_{\theta }}}{\partial {{r}^{2}}}+\frac{2}{{{r}^{2}}}\frac{{{\partial }^{2}}{{u}_{\theta }}}{\partial {{\theta }^{2}}}+\frac{1}{r}\frac{\partial {{u}_{\theta }}}{\partial r}-\frac{{{u}_{\theta }}}{{{r}^{2}}}-\frac{\partial {{P}_{r\theta }}}{\partial r}-\frac{\partial {{P}_{\theta r}}}{\partial r}-\frac{2}{r}\frac{\partial {{P}_{\theta \theta }}}{\partial \theta }-2\frac{{{P}_{r\theta }}+{{P}_{\theta r}}}{r} \right)\\
\nonumber
 &+{{\lambda }_{\rm e}}\left( \frac{1}{r}\frac{{{\partial }^{2}}{{u}_{r}}}{\partial r\partial \theta }+\frac{1}{{{r}^{2}}}\frac{{{\partial }^{2}}{{u}_{\theta }}}{\partial {{\theta }^{2}}}+\frac{1}{{{r}^{2}}}\frac{\partial {{u}_{r}}}{\partial \theta }-\frac{1}{r}\frac{\partial {{P}_{rr}}}{\partial \theta }-\frac{1}{r}\frac{\partial {{P}_{\theta \theta }}}{\partial \theta } \right)\\
 &+{{\mu }_{\rm c}}\left( \frac{1}{{{r}^{2}}}\frac{\partial {{u}_{r}}}{\partial \theta }-\frac{{{u}_{\theta }}}{{{r}^{2}}}-\frac{1}{r}\frac{{{\partial }^{2}}{{u}_{r}}}{\partial r\partial \theta }+\frac{1}{r}\frac{\partial {{u}_{\theta }}}{\partial r}+\frac{{{\partial }^{2}}{{u}_{\theta }}}{\partial {{r}^{2}}}-\frac{\partial {{P}_{\theta r}}}{\partial r}+\frac{\partial {{P}_{r\theta }}}{\partial r} \right).
\end{align}
}
In addition, the Curl of the micro-distortion tensor can be written as \cite{reddy2013introduction}
{\small
\begin{align}
\text{Curl} \, P=\left( \frac{\partial {{P}_{\theta r}}}{\partial r}-\frac{1}{r}\frac{\partial {{P}_{rr}}}{\partial \theta }+\frac{1}{r}({{P}_{r\theta }}+{{P}_{\theta r}}) \right)\,{{\mathbf{e}}_{\mathbf{z}}}\otimes{{\mathbf{e}}_{\mathbf{r}}}+\left( \frac{\partial {{P}_{\theta \theta }}}{\partial r}+\frac{{{P}_{\theta \theta }}-{{P}_{rr}}}{r}-\frac{1}{r}\frac{\partial {{P}_{r\theta }}}{\partial \theta } \right)\,{{\mathbf{e}}_{\mathbf{z}}}\otimes{{\mathbf{e}}_{\mathbf{\theta }}} \,, 
\label{Curl_P} 
\end{align}
}
in which ${\mathbf{e}}_{\mathbf{z}}$ denotes the unit vector normal to the $r-\theta$ plane. As a result, using Eqs. (\ref{constEq}) and (\ref{Curl_P}), the nonsymmetric moment tensor $m$ is obtained as

{\small
\begin{align}
   m&={{m}_{zr}}\,{{\mathbf{e}}_{\mathbf{z}}}\otimes{{\mathbf{e}}_{\mathbf{r}}}+{{m}_{z\theta }}\,{{\mathbf{e}}_{\mathbf{z}}}\otimes{{\mathbf{e}}_{\mathbf{\theta }}} \\  \nonumber
 & ={{\mu }_{\textrm{macro}}}L_{c}^{2}\left( \frac{\partial {{P}_{\theta r}}}{\partial r}-\frac{1}{r}\frac{\partial {{P}_{rr}}}{\partial \theta }+\frac{1}{r}({{P}_{r\theta }}+{{P}_{\theta r}}) \right)\,{{\mathbf{e}}_{\mathbf{z}}}\otimes{{\mathbf{e}}_{\mathbf{r}}}+{{\mu }_{\textrm{macro}}}L_{c}^{2}\left( \frac{\partial {{P}_{\theta \theta }}}{\partial r}+\frac{{{P}_{\theta \theta }}-{{P}_{rr}}}{r}-\frac{1}{r}\frac{\partial {{P}_{r\theta }}}{\partial \theta } \right)\,{{\mathbf{e}}_{\mathbf{z}}}\otimes{{\mathbf{e}}_{\mathbf{\theta }}} \, . 
 \label{m_mat} 
 \end{align}
 }
 Furthermore, the Curl of $m$ is derived as
{\small
\begin{align}
\nonumber \frac{1}{{{\mu }_{\textrm{macro}}}L_{c}^{2}}{\text{Curl }m}=&\frac{1}{{{\mu }_{\textrm{macro}}}L_{c}^{2}}[ \left( \frac{1}{r}\frac{\partial {{m}_{zr}}}{\partial \theta }-\frac{{{m}_{z\theta }}}{r} \right){{\mathbf{e}}_{\mathbf{r}}}\otimes{{\mathbf{e}}_{\mathbf{r}}}-\left( \frac{\partial {{m}_{z\theta }}}{\partial r} \right){{\mathbf{e}}_{\mathbf{\theta }}}\otimes{{\mathbf{e}}_{\mathbf{\theta }}} \\ \nonumber
&\quad \quad \quad \quad +\left( \frac{1}{r}\frac{\partial {{m}_{z\theta }}}{\partial \theta }+\frac{{{m}_{zr}}}{r} \right){{\mathbf{e}}_{\mathbf{r}}}\otimes{{\mathbf{e}}_{\mathbf{\theta }}}-\left( \frac{\partial {{m}_{zr}}}{\partial r} \right){{\mathbf{e}}_{\mathbf{\theta }}}\otimes{{\mathbf{e}}_{\mathbf{r}}}] \\ 
  =&\left( \frac{1}{r}\frac{{{\partial }^{2}}{{P}_{\theta r}}}{\partial r\partial \theta }-\frac{1}{{{r}^{2}}}\frac{{{\partial }^{2}}{{P}_{rr}}}{\partial {{\theta }^{2}}}+\frac{1}{{{r}^{2}}}\frac{\partial {{P}_{\theta r}}}{\partial \theta }-\frac{1}{r}\frac{\partial {{P}_{\theta \theta }}}{\partial r}-\frac{{{P}_{\theta \theta }}-{{P}_{rr}}}{{{r}^{2}}}+\frac{2}{{{r}^{2}}}\frac{\partial {{P}_{r\theta }}}{\partial \theta } \right){{\mathbf{e}}_{\mathbf{r}}}\otimes{{\mathbf{e}}_{\mathbf{r}}} \\  \nonumber
 & +\left( \frac{1}{r}\frac{{{\partial }^{2}}{{P}_{\theta \theta }}}{\partial r\partial \theta }+\frac{1}{{{r}^{2}}}\frac{\partial {{P}_{\theta \theta }}}{\partial \theta }-\frac{2}{{{r}^{2}}}\frac{\partial {{P}_{rr}}}{\partial \theta }-\frac{1}{{{r}^{2}}}\frac{{{\partial }^{2}}{{P}_{r\theta }}}{\partial {{\theta }^{2}}}+\frac{1}{r}\frac{\partial {{P}_{\theta r}}}{\partial r}+\frac{{{P}_{r\theta }}+{{P}_{\theta r}}}{{{r}^{2}}} \right){{\mathbf{e}}_{\mathbf{r}}}\otimes{{\mathbf{e}}_{\mathbf{\theta }}} \\ \nonumber
 & +\left( -\frac{{{\partial }^{2}}{{P}_{\theta r}}}{\partial {{r}^{2}}}-\frac{1}{{{r}^{2}}}\frac{\partial {{P}_{rr}}}{\partial \theta }+\frac{1}{r}\frac{{{\partial }^{2}}{{P}_{rr}}}{\partial r\partial \theta }+\frac{{{P}_{r\theta }}+{{P}_{\theta r}}}{{{r}^{2}}}-\frac{1}{r}\frac{\partial {{P}_{r\theta }}}{\partial r}-\frac{1}{r}\frac{\partial {{P}_{\theta r}}}{\partial r} \right){{\mathbf{e}}_{\mathbf{\theta }}}\otimes{{\mathbf{e}}_{\mathbf{r}}} \\  
 & +\left( -\frac{{{\partial }^{2}}{{P}_{\theta \theta }}}{\partial {{r}^{^{2}}}}+\frac{{{P}_{\theta \theta }}-{{P}_{rr}}}{{{r}^{2}}}-\frac{1}{r}\frac{\partial {{P}_{\theta \theta }}}{\partial r}+\frac{1}{r}\frac{\partial {{P}_{rr}}}{\partial r}-\frac{1}{{{r}^{2}}}\frac{\partial {{P}_{r\theta }}}{\partial \theta }+\frac{1}{r}\frac{{{\partial }^{2}}{{P}_{r\theta }}}{\partial r\partial \theta } \right){{\mathbf{e}}_{\mathbf{\theta }}}\otimes{{\mathbf{e}}_{\mathbf{\theta }}} \,. \nonumber
\end{align}
}

Substituting Eqs. (\ref{constEq}), (\ref{eq:P_polar}), (\ref{grad_u}), and (15) in Eq. (\ref{eq:equi_RM}), yields
{\small
\begin{align}
\nonumber \eta \frac{{{\partial }^{2}}{{P}_{rr}}}{\partial {{t}^{2}}}=2{{\mu }_{\rm e}}&\left( \frac{\partial {{u}_{r}}}{\partial r}-{{P}_{rr}} \right)+{{\lambda }_{\rm e}}\left( \frac{\partial {{u}_{r}}}{\partial r}+\frac{1}{r}\left( \frac{\partial {{u}_{\theta }}}{\partial \theta }+{{u}_{r}} \right)-({{P}_{rr}}+{{P}_{\theta \theta }}) \right) \\  \nonumber
 & -2{{\mu }_{\rm micro}}{{P}_{rr}}-{{\lambda }_{\rm micro}}({{P}_{rr}}+{{P}_{\theta \theta }}) \\ 
 &-{{\mu }_{\textrm{macro}}}L_{c}^{2}\left( \frac{1}{r}\frac{{{\partial }^{2}}{{P}_{\theta r}}}{\partial r\partial \theta }-\frac{1}{{{r}^{2}}}\frac{{{\partial }^{2}}{{P}_{rr}}}{\partial {{\theta }^{2}}}+\frac{1}{{{r}^{2}}}\frac{\partial {{P}_{\theta r}}}{\partial \theta }-\frac{1}{r}\frac{\partial {{P}_{\theta \theta }}}{\partial r}-\frac{{{P}_{\theta \theta }}-{{P}_{rr}}}{{{r}^{2}}}+\frac{2}{{{r}^{2}}}\frac{\partial {{P}_{r\theta }}}{\partial \theta } \right),
 \label{govP1} 
\end{align}

\begin{align}
\nonumber \eta \frac{{{\partial }^{2}}{{P}_{r\theta }}}{\partial {{t}^{2}}} =2{{\mu }_{\rm e}}&\left( \frac{1}{2r}\left( \frac{\partial {{u}_{r}}}{\partial \theta }-{{u}_{\theta }} \right)+\frac{1}{2}\frac{\partial {{u}_{\theta }}}{\partial r}-\frac{{{P}_{r\theta }}+{{P}_{\theta r}}}{2} \right) \\  \nonumber
 & -{{\mu }_{\rm micro}}({{P}_{r\theta }}+{{P}_{\theta r}})+2{{\mu }_{\rm c}}\left( \frac{1}{2r}\left( \frac{\partial {{u}_{r}}}{\partial \theta }-{{u}_{\theta }} \right)-\frac{1}{2}\frac{\partial {{u}_{\theta }}}{\partial r}-\frac{{{P}_{r\theta }}-{{P}_{\theta r}}}{2} \right) \\  
  & -{{\mu }_{\textrm{macro}}}L_{c}^{2}\left( \frac{1}{r}\frac{{{\partial }^{2}}{{P}_{\theta \theta }}}{\partial r\partial \theta }+\frac{1}{{{r}^{2}}}\frac{\partial {{P}_{\theta \theta }}}{\partial \theta }-\frac{2}{{{r}^{2}}}\frac{\partial {{P}_{rr}}}{\partial \theta }-\frac{1}{{{r}^{2}}}\frac{{{\partial }^{2}}{{P}_{r\theta }}}{\partial {{\theta }^{2}}}+\frac{1}{r}\frac{\partial {{P}_{\theta r}}}{\partial r}+\frac{{{P}_{r\theta }}+{{P}_{\theta r}}}{{{r}^{2}}} \right),
  \label{govP2} 
\end{align}

\begin{align}
\nonumber  \eta \frac{{{\partial }^{2}}{{P}_{\theta r}}}{\partial {{t}^{2}}}=2{{\mu }_{\rm e}}&\left( \frac{1}{2r}\left( \frac{\partial {{u}_{r}}}{\partial \theta }-{{u}_{\theta }} \right)+\frac{1}{2}\frac{\partial {{u}_{\theta }}}{\partial r}-\frac{{{P}_{r\theta }}+{{P}_{\theta r}}}{2} \right) \\ \nonumber
 & -{{\mu }_{\rm micro}}({{P}_{r\theta }}+{{P}_{\theta r}})+2{{\mu }_{c}}\left( -\frac{1}{2r}\left( \frac{\partial {{u}_{r}}}{\partial \theta }-{{u}_{\theta }} \right)+\frac{1}{2}\frac{\partial {{u}_{\theta }}}{\partial r}-\frac{{{P}_{\theta r}}-{{P}_{r\theta }}}{2} \right) \\  
 & -{{\mu }_{\textrm{macro}}}L_{c}^{2}\left( -\frac{{{\partial }^{2}}{{P}_{\theta r}}}{\partial {{r}^{2}}}-\frac{1}{{{r}^{2}}}\frac{\partial {{P}_{rr}}}{\partial \theta }+\frac{1}{r}\frac{{{\partial }^{2}}{{P}_{rr}}}{\partial r\partial \theta }+\frac{{{P}_{r\theta }}+{{P}_{\theta r}}}{{{r}^{2}}}-\frac{1}{r}\frac{\partial {{P}_{r\theta }}}{\partial r}-\frac{1}{r}\frac{\partial {{P}_{\theta r}}}{\partial r} \right),
  \label{govP3} 
\end{align} 

\begin{align}
 \nonumber \eta \frac{{{\partial }^{2}}{{P}_{\theta \theta }}}{\partial {{t}^{2}}} =2{{\mu }_{\rm e}}&\left( \frac{1}{r}\left( \frac{\partial {{u}_{\theta }}}{\partial \theta }+{{u}_{r}} \right)-{{P}_{\theta \theta }} \right)+{{\lambda }_{\rm e}}\left( \frac{\partial {{u}_{r}}}{\partial r}+\frac{1}{r}\left( \frac{\partial {{u}_{\theta }}}{\partial \theta }+{{u}_{r}} \right)-({{P}_{rr}}+{{P}_{\theta \theta }}) \right) \\ \nonumber
 & -2{{\mu }_{\rm micro}}{{P}_{\theta \theta }}-{{\lambda }_{\rm micro}}({{P}_{rr}}+{{P}_{\theta \theta }}) \\
 &-{{\mu }_{\textrm{macro}}}L_{c}^{2}\left( -\frac{{{\partial }^{2}}{{P}_{\theta \theta }}}{\partial {{r}^{^{2}}}}+\frac{{{P}_{\theta \theta }}-{{P}_{rr}}}{{{r}^{2}}}-\frac{1}{r}\frac{\partial {{P}_{\theta \theta }}}{\partial r}+\frac{1}{r}\frac{\partial {{P}_{rr}}}{\partial r}-\frac{1}{{{r}^{2}}}\frac{\partial {{P}_{r\theta }}}{\partial \theta }+\frac{1}{r}\frac{{{\partial }^{2}}{{P}_{r\theta }}}{\partial r\partial \theta } \right). 
 \label{govP4} 
\end{align}
}

\subsection{Consistent coupling boundary conditions in polar coordinates}  
Boundary conditions are an essential part of any model formulation. For classical linear elasticity, we typically have either Dirichlet and  Neumann conditions on the displacement or stress, respectively. One of the challenges in generalized continua is the inclusion of additional degrees of freedom. For them, the interpretation and adequacy of boundary conditions is an ongoing field of research \cite{d2022consistent}. From the physical perspective, if one prefers to understand the microdistortion $P$, it is considered necessary to prescribe the full Dirichlet condition at the boundary. However, in the relaxed micromorphic model, this is not possible because only the curl of $P$ is controlled within the energy function. Therefore, one only needs to prescribe the tangential traces of $P$ at the boundary. However, it is not clear at all if $P \cdot \tau$ should be set to zero. That the tangential traces of $P$ need to be prescribed is evident from the mathematical existence theories based on the incompatible Korn's inequality \cite{lewintan2021korn,neff2015poincare,gmeineder2023optimal,gmeineder2024korn}. Neff and co-workers \cite{d2022consistent} have observed that at the Dirichlet part of the boundary, where the displacement $u$ is prescribed, there is also an additional prescribed tangential part of $\text{D}u \cdot \tau$. The idea of the consistent coupling boundary condition is to set $P \cdot \tau=\text{D}u \cdot \tau$ on the boundary. Since $\text{Curl}(\text{D}u - P) = -\text{Curl}(P)$, the combination $(\text{D}u - P)$ indeed has tangential traces at the boundary. 

Here, we rewrite the consistent coupling condition in polar coordinates for the general case. In this context, $\tau = \tau_r \,\mathbf{e}_{\mathbf{r}} + \tau_{\theta} \,\mathbf{e}_\theta$, and applying Eqs. (\ref{eq:P_polar}) and (\ref{grad_u}) yields
{\small
\begin{align}
P \cdot\tau =({{P}_{rr}}{{\tau }_{r}}+{{P}_{r\theta }}{{\tau }_{\theta }})\,{{\mathbf{e}}_{\mathbf{r}}}+({{P}_{\theta r}}{{\tau }_{r}}+{{P}_{\theta \theta }}{{\tau }_{\theta }})\,{{\mathbf{e}}_{\theta}},
\label{CCBC1} 
\end{align}

\begin{align}
\text{D}u\cdot \tau =(\frac{\partial {{u}_{r}}}{\partial r}{{\tau }_{r}}+\frac{1}{r}\left( \frac{\partial {{u}_{r}}}{\partial \theta }-{{u}_{\theta }} \right){{\tau }_{\theta }})\,{{\mathbf{e}}_{\mathbf{r}}}+(\frac{\partial {{u}_{\theta }}}{\partial r}{{\tau }_{r}}+\frac{1}{r}\left( \frac{\partial {{u}_{\theta }}}{\partial \theta }+{{u}_{r}} \right)\,{{\tau }_{\theta }})\,{{\mathbf{e}}_{\theta}} .
\label{CCBC2} 
\end{align}
}
By substituting Eqs. (\ref{CCBC1}) and (\ref{CCBC2}) into $P \cdot \tau=\text{D}u \cdot \tau$ at the boundary of the domain, we derive the following result:
{\small
\begin{align}
{{P}_{rr}}{{\tau }_{r}}+{{P}_{r\theta }}{{\tau }_{\theta }}=\frac{\partial {{u}_{r}}}{\partial r}{{\tau }_{r}}+\frac{1}{r}\left( \frac{\partial {{u}_{r}}}{\partial \theta }-{{u}_{\theta }} \right){{\tau }_{\theta }}, \\ {{P}_{\theta r}}{{\tau }_{r}}+{{P}_{\theta \theta }}{{\tau }_{\theta }}=\frac{\partial {{u}_{\theta }}}{\partial r}{{\tau }_{r}}+\frac{1}{r}\left( \frac{\partial {{u}_{\theta }}}{\partial \theta }+{{u}_{r}} \right){{\tau }_{\theta }} .
\label{CCBpol1} 
\end{align}
}
If the boundary is circular with radius $R$, we have $\tau = \mathbf{e}_\theta$. Therefore, $\tau_{r} = 0$ and $\tau_{\theta} = 1$, leading to the following result:
{\small
\begin{align}
P \cdot\tau ={{P}_{r\theta }}\,{{\mathbf{e}}_{\mathbf{r}}}+{{P}_{\theta \theta }}\,{{\mathbf{e}}_{\theta}},
\label{CCBC1-cir} 
\end{align}

\begin{align}
\text{D}u\cdot \tau =\frac{1}{r}\left( \frac{\partial {{u}_{r}}}{\partial \theta }-{{u}_{\theta }} \right)\,{{\mathbf{e}}_{\mathbf{r}}}+\frac{1}{r}\left( \frac{\partial {{u}_{\theta }}}{\partial \theta }+{{u}_{r}} \right)\,{{\mathbf{e}}_{\theta}} .
\label{CCBC2-cir} 
\end{align}
}
Equating Eqs. (\ref{CCBC1-cir}) and (\ref{CCBC2-cir}), we obtain 
{\small
\begin{align}
{{P}_{r\theta }}(R)=\frac{1}{R}\left( \frac{\partial {{u}_{r}(R)}}{\partial \theta }-{{u}_{\theta }(R)} \right),
\label{CCBpol-cir1} 
\end{align}

\begin{align}
{{P}_{\theta \theta }}(R)=\frac{1}{R}\left( \frac{\partial {{u}_{\theta }(R)}}{\partial \theta }+{{u}_{r}(R)} \right).
\label{CCBpol-cir2} 
\end{align}
}

\section{Elastostatic axisymmetric extension problem}
\subsection{Formulations and solution method}
Here, we consider an elastostatic axisymmetric extension problem that involves analyzing the deformation of a long circular cylinder with radius $R$ under the influence of uniform radial displacement at the boundary (see Figure \ref{figure1}). In this special case, both the displacement field and the micro-distortion tensor depend solely on the radial coordinate ($u_r(r)$, $P_{rr}(r)$, $P_{\theta \theta}(r)$, $P_{r \theta}(r)$ and $P_{\theta r}(r)$), and we have also ${{u}_{\theta }}=0$. By dropping the time-dependent terms and using the mentioned assumptions, the governing equations are simplified as follows: 
{\small
\begin{align}
  (2{{\mu }_{\rm e}}+{{\lambda }_{\rm e}})\frac{d}{dr}\left( \frac{d{{u}_{r}}}{dr}-{{P}_{rr}} \right)+{{\lambda }_{\rm e}}\frac{d}{dr}\left( \frac{{{u}_{r}}}{r}-{{P}_{\theta \theta }} \right)=-2{{\mu }_{\rm e}}\frac{d}{dr}\left( \frac{{{u}_{r}}}{r}-{{P}_{\theta \theta }} \right)-2{{\mu }_{\rm e}}\left( \frac{d{{P}_{\theta \theta }}}{dr}+\frac{{{P}_{\theta \theta }}-{{P}_{rr}}}{r} \right) ,
  \label{EAEP1}
\end{align}

\begin{align}
{{\mu }_{\text{e}}}\left( \frac{d{{P}_{r\theta }}}{dr}+\frac{d{{P}_{\theta r}}}{dr}+2\frac{{{P}_{r\theta }}+{{P}_{\theta r}}}{r} \right)={{\mu }_{\text{c}}}\left( \frac{d{{P}_{r\theta }}}{dr}-\frac{d{{P}_{\theta r}}}{dr} \right),
  \label{EAEP2}
\end{align}

\begin{align}
(2{{\mu }_{\rm e}}+{{\lambda }_{\rm e}})\left( \frac{d{{u}_{r}}}{dr}-{{P}_{rr}} \right)+{{\lambda }_{\rm e}}\left( \frac{{{u}_{r}}}{r}-{{P}_{\theta \theta }} \right)=2{{\mu }_{\rm micro}}{{P}_{rr}}+{{\lambda }_{\rm micro}}({{P}_{rr}}+{{P}_{\theta \theta }})-\frac{{{\mu }_{\textrm{macro}}}L_{c}^{2}}{r}\left( \frac{d{{P}_{\theta \theta }}}{dr}+\frac{{{P}_{\theta \theta }}-{{P}_{rr}}}{r} \right) ,
\label{EAEP3}
\end{align}

\begin{align}
{{\mu }_{\text{e}}}\left( {{P}_{r\theta }}+{{P}_{\theta r}} \right)+{{\mu }_{\text{micro}}}({{P}_{r\theta }}+{{P}_{\theta r}})+{{\mu }_{\text{c}}}\left( {{P}_{r\theta }}-{{P}_{\theta r}} \right)+{{\mu }_{\text{macro}}}L_{\text{c}}^{\text{2}}\left( \frac{1}{r}\frac{d{{P}_{\theta r}}}{dr}+\frac{{{P}_{r\theta }}+{{P}_{\theta r}}}{{{r}^{2}}} \right)=0,
\label{EAEP4}
\end{align}

\begin{align}
{{\mu }_{\text{e}}}\left( {{P}_{r\theta }}+{{P}_{\theta r}} \right)+{{\mu }_{\text{micro}}}({{P}_{r\theta }}+{{P}_{\theta r}})+{{\mu }_{\text{c}}}\left( {{P}_{\theta r}}-{{P}_{r\theta }} \right)+{{\mu }_{\text{macro}}}L_{\text{c}}^{\text{2}}\left( -\frac{{{d}^{2}}{{P}_{\theta r}}}{d{{r}^{2}}}+\frac{{{P}_{r\theta }}+{{P}_{\theta r}}}{{{r}^{2}}}-\frac{1}{r}\frac{d{{P}_{r\theta }}}{dr}-\frac{1}{r}\frac{d{{P}_{\theta r}}}{dr} \right)=0,
\label{EAEP5}
\end{align}

\begin{align}
(2{{\mu }_{\rm e}}+{{\lambda }_{\rm e}})\left( \frac{{{u}_{r}}}{r}-{{P}_{\theta \theta }} \right)+{{\lambda }_{\rm e}}\left( \frac{d{{u}_{r}}}{dr}-{{P}_{rr}} \right)=2{{\mu }_{\rm micro}}{{P}_{\theta \theta }}+{{\lambda }_{\rm micro}}({{P}_{rr}}+{{P}_{\theta \theta }})-{{\mu }_{\textrm{macro}}}L_{c}^{2}\frac{d}{dr}\left( \frac{d{{P}_{\theta \theta }}}{dr}+\frac{{{P}_{\theta \theta }}-{{P}_{rr}}}{r} \right)  .
\label{EAEP6}
\end{align}
}

\begin{figure}[!h]
\centering
\includegraphics[scale=0.2]{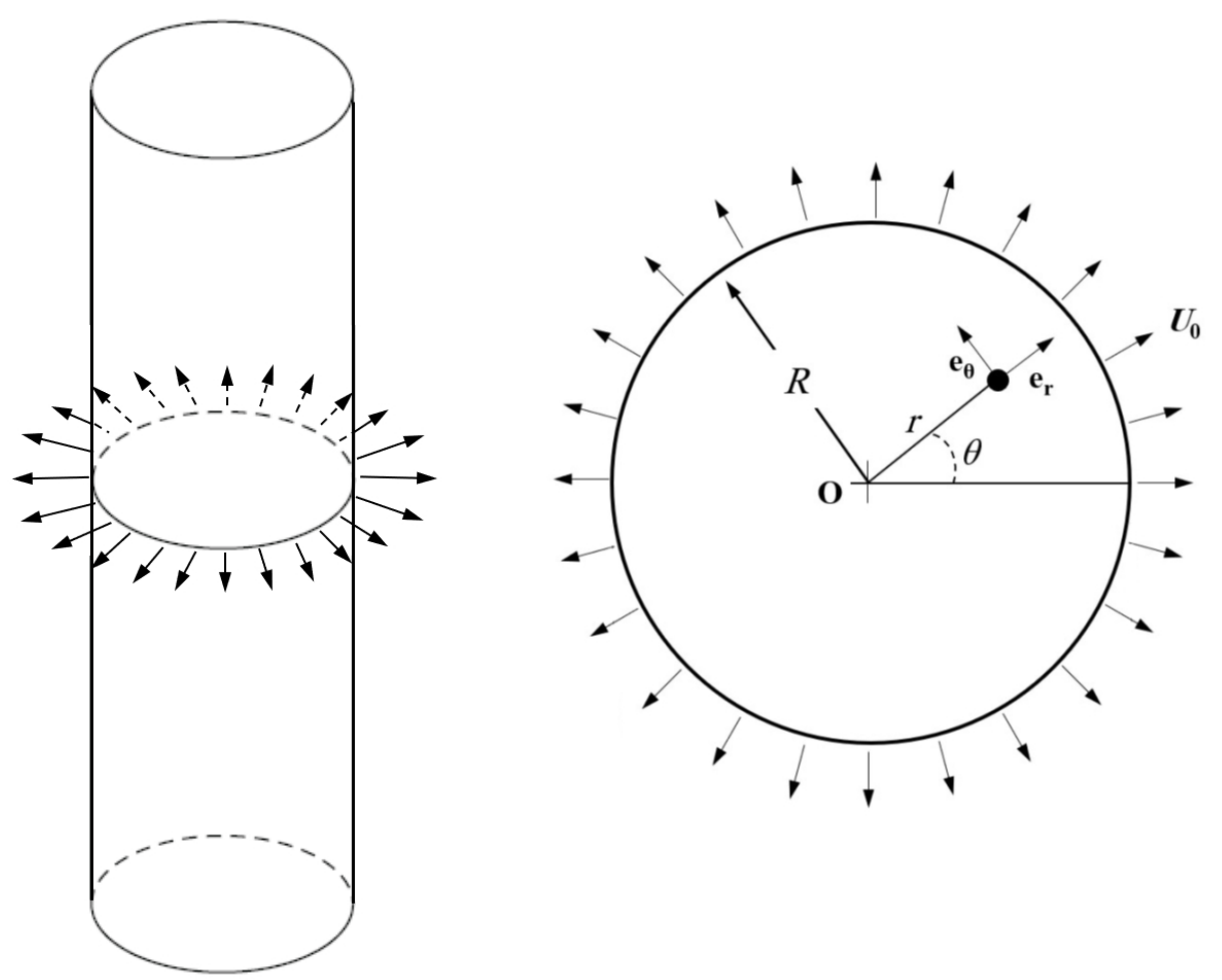}
\caption{Cross-section of the long circular cylinder with radius $R$ under the influence of uniform radial displacement $U_0$. Due to the length of the cylinder, plane-strain conditions are warranted.}
\label{figure1}
\end{figure}

The Dirichlet boundary conditions for the axisymmetric extension is
{\small
\begin{align}
{{u}_{r}}(R)={{U}_{0}}.
\label{boudcon1}
\end{align}
}
Furthermore, using Eqs. (\ref{CCBpol-cir1}) and (\ref{CCBpol-cir2}), the consistent coupling boundary conditions for the axisymmetric problem are as follows: 
{\small
\begin{align}
{{P}_{r\theta }}(R)=0, \quad \quad {{P}_{\theta \theta }}(R)=\frac{{{U}_{0}}}{R}.
\label{boudcon2}
\end{align}
}
To enhance readability, in the following, we utilize the shortened forms: $\mu_{\rm M} \equiv \mu_{\rm macro}$, $\mu_{\rm m} \equiv \mu_{\rm micro}$,  $\lambda_{\textrm{m}} \equiv \lambda_{\textrm{micro}}$ and $\lambda_{\textrm{M}} \equiv \lambda_{\textrm{macro}}$. Since the cylinder is long, the plane-strain conditions are warranted. In the context of plane-strain conditions, the bulk micro-moduli  $\kappa_{\rm e}$ and $\kappa_{\rm m} \equiv\kappa_{\rm micro}$ are connected to their corresponding Lamé-type micro-moduli via the two-dimensional relationships \cite{gourgiotis2024green}:
{\small
\begin{equation}\label{eq:KAPPA}
\kappa_{\rm e}\coloneqq \lambda_{\rm e}+\mu_{\rm e}, \qquad \quad \kappa_{\rm m} \coloneqq\lambda_{\rm m}+\mu_{\rm m}\,.
\end{equation}
}
Accordingly, the relations between the macro moduli ($\mu_{\rm M},\lambda_{\rm M}, \kappa_{\rm M}$) and the micro-moduli in the plane strain case become \cite{neff2007geometrically}
{\small
\begin{align}
\label{defmod}
\mu_{\rm M} &\coloneqq \dfrac{\mu_{\rm e} \, \mu_{\rm m}}{\mu_{\rm e} + \mu_{\rm m}} 
\qquad\Leftrightarrow\qquad \frac{1}{\mu_{\rm M}}=\frac{1}{\mu_{\rm e}}+\frac{1}{\mu_{\rm m}}\,,\notag\\
\kappa_{\rm M} &\coloneqq \dfrac{\kappa_{\rm e} \, \kappa_{\rm m}}{\kappa_{\rm e} + \kappa_{\rm m}} \qquad\Leftrightarrow\qquad
 \frac{1}{\kappa_{\rm M}}=\frac{1}{\kappa_{\rm e}}+\frac{1}{\kappa_{\rm m}}\,,\\
\lambda_{\rm m} &\coloneqq
\dfrac{(\mu_{\rm e} + \lambda_{\rm e})(\mu_{\rm m} + \lambda_{\rm m})}{(\mu_{\rm e} + \lambda_{\rm e}) + (\mu_{\rm m} + \lambda_{\rm m})} -\dfrac{\mu_{\rm e} \, \mu_{\rm m}}{\mu_{\rm e} + \mu_{\rm m}}
\, ,
\notag
\end{align}
}
where $\kappa_{\rm M} \equiv \kappa_{\rm macro}$ with $\kappa_{\rm M}= \mu_{\rm M}+\lambda_{\rm M} $. To obtain the analytical solution for the problem, it is essential to reformulate the system (\ref{EAEP1})-(\ref{EAEP6}). First, Eq. (\ref{EAEP1}) is rearranged as follows:
{\small
\begin{align}
  & (2{{\mu }_{\rm e}}+{{\lambda }_{\rm e}})\frac{d}{dr}\left( \frac{d{{u}_{r}}}{dr}+\frac{{{u}_{r}}}{r} \right)-2{{\mu }_{\rm e}}(\frac{d{{P}_{rr}}}{dr}-\frac{{{P}_{\theta \theta }}-{{P}_{rr}}}{r})-{{\lambda }_{\rm e}}\frac{d}{dr}({{P}_{\theta \theta }}+{{P}_{rr}})=0 .
  \label{EAEP_ref}
\end{align}
}
Then, by differentiating Eq. (\ref{EAEP3}) with respect to $r$, and using Eqs. (\ref{EAEP1}) and (\ref{EAEP6}) along with the simplification, the resultant equation is obtained as:
{\small
\begin{align}
2{{\mu }_{\rm m}}\left( \frac{d{{P}_{rr}}}{dr}+\frac{{{P}_{rr}}-{{P}_{\theta \theta }}}{r} \right)+{{\lambda }_{\rm m}}\frac{d}{dr}({{P}_{\theta \theta }}+{{P}_{rr}})=0 .
\label{ref_EAEP2}
\end{align}
}
Substituting Eq. (\ref{ref_EAEP2}) into Eq. (\ref{EAEP_ref}) yields
{\small
\begin{align}
\frac{d{{u}_{r}}}{dr}+\frac{{{u}_{r}}}{r}=\frac{{{\mu }_{\rm m}}{{\lambda }_{\rm e}}-{{\mu }_{\rm e}}{{\lambda }_{\rm m}}}{{{\mu }_{\rm m}}({{\mu }_{\rm e}}+{{\kappa}_{\rm e}})}({{P}_{\theta \theta }}+{{P}_{rr}})+{{C}_{1}} ,
\label{Int_EAEP}
\end{align}
}
where $C_1$ is an integration constant. In addition, by summing Eqs. (\ref{EAEP3}) and (\ref{EAEP6}) we have
{\small
\begin{align}
2({{\mu }_{\rm e}}+{{\lambda }_{\rm e}})\left( \frac{d{{u}_{r}}}{dr}+\frac{{{u}_{r}}}{r} \right)=2({\kappa}_{\rm m}+{\kappa}_{\rm e})({{P}_{rr}}+{{P}_{\theta \theta }})-{{\mu }_{\rm M}}L_{c}^{2}\frac{d}{dr}\left( \frac{d{{P}_{\theta \theta }}}{dr}+\frac{{{P}_{\theta \theta }}-{{P}_{rr}}}{r} \right)-{{\mu }_{\rm M}}L_{c}^{2}\frac{1}{r}\left( \frac{d{{P}_{\theta \theta }}}{dr}+\frac{{{P}_{\theta \theta }}-{{P}_{rr}}}{r} \right).
\label{EAEP3_ref}
\end{align}
}
Now, we define the following variables:
{\small
\begin{align}
X=\frac{d{{u}_{r}}}{dr}-{{P}_{rr}}, \quad \quad Y=\frac{{{u}_{r}}}{r}-{{P}_{\theta \theta }}, \quad \quad Z={{P}_{\theta \theta }}+{{P}_{rr}} .
\label{new_var}
\end{align}
}
Substituting Eq. (\ref{new_var}) into Eqs. (\ref{ref_EAEP2})-(\ref{EAEP3_ref}), we obtain
{\small
\begin{align}
\frac{dY}{dr}+\frac{Y-X}{r}&=-\left(\frac{{{\kappa}_{\rm m}}+{{\mu }_{\rm m}}}{2{{\mu }_{\rm m}}}\right)\frac{dZ}{dr}, \\
X+Y&=-\frac{{{\mu }_{\rm e}}({{\kappa}_{\rm m}}+{{\mu }_{\rm m}})}{{{\mu }_{\rm m}}({{\kappa}_{\rm e}}+{{\mu }_{\rm e}})}Z+{{C}_{1}},\\
X+Y&=\frac{{{\kappa}_{\rm m}}}{{{\kappa}_{\rm e}}}Z+\frac{{{\mu }_{\rm M}}L_{c}^{2}}{2{{\kappa}_{\rm e}}}\frac{d}{dr}\left( \frac{dY}{dr}+\frac{Y-X}{r}\right)+\frac{{{\mu }_{\rm M}}L_{c}^{2}}{2{{\kappa}_{\rm e}}}\frac{1}{r}\left( \frac{dY}{dr}+\frac{Y-X}{r} \right).
\end{align}
}
By using Eqs. (44) and (45) and simplifying the equations, Eq. (46) can be rewritten as
{\small
\begin{align}
\frac{{{d}^{2}}Z}{d{{r}^{2}}}+\frac{1}{r}\frac{dZ}{dr}-aZ+b=0,
\label{EQZ}
\end{align}
}
where
{\small
\begin{align}
a=\frac{4}{{{\mu }_{\rm M}}L_{c}^{2}}\left( \frac{{{\mu }_{\rm e}}{{\kappa}_{\rm e}}}{{{\kappa}_{\rm e}}+{{\mu}_{\rm e}}}+\frac{{{\mu }_{\rm m}}{{\kappa}_{\rm m}}}{{{\kappa}_{\rm m}}+{{\mu}_{\rm m}}} \right), \quad \quad b=\frac{4}{{{\mu }_{\rm M}}L_{c}^{2}}\frac{{{C}_{1}}{{\kappa}_{\rm e}}{{\mu }_{\rm m}}}{{{\kappa}_{\rm m}}+{{\mu }_{\rm m}}}.
\label{equa}
\end{align}
}
The general solution of Eq. (\ref{EQZ}) is given by
{\small
\begin{align}
Z(r)=\frac{b}{a}+{{D}_{1}}{{I}_{0}}(\sqrt{a}r)+{{D}_{2}}{{K}_{0}}(\sqrt{a}r),
\label{Sol}
\end{align}
}
where $D_1$ and $D_2$ are unknown constants, and $I_0$ and $K_0$ are modified Bessel functions of the first and second kind, respectively, of order zero. To ensure the solution remains finite at the center of the plate, $D_2$ must be set to 0 to eliminate the infinite value of $K_0$ when $r = 0$. As a result, the solution becomes
{\small
\begin{align}
Z(r)=\frac{b}{a}+{{D}_{1}}{{I}_{0}}(\sqrt{a}r).
\label{NewSol}
\end{align}
}
Using Eq. (\ref{new_var}) and (\ref{NewSol}), Eq. (\ref{Int_EAEP}) can be rewritten as
{\small
\begin{align}
\frac{d{{u}_{r}}}{dr}+\frac{{{u}_{r}}}{r}={{C}_{1}}A+{{D}_{1}}B{{I}_{0}}(\sqrt{a}r), 
\label{EQu}
\end{align}
}
where now
{\small
\begin{align}
A=1+\frac{({{\mu }_{\rm m}}{{\lambda }_{\rm e}}-{{\mu }_{\rm e}}{{\lambda }_{\rm m}}){{\kappa}_{\rm e}}}{{{\mu }_{\rm e}}({{\kappa}_{\rm m}}+{{\mu }_{\rm m}}){{\kappa}_{\rm e}}+{{\mu }_{\rm m}}({{\kappa}_{\rm e}}+{{\mu}_{\rm e}}){{\kappa}_{\rm m}}}, \quad \quad B=\frac{{{\mu }_{\rm m}}{{\lambda }_{\rm e}}-{{\mu }_{\rm e}}{{\lambda }_{\rm m}}}{{{\mu }_{\rm m}}({{\mu }_{\rm e}}+{{\kappa }_{\rm e}})}.
\label{A_B}
\end{align}
}

The general solution of Eq. (\ref{EQu}) is obtained by
{\small
\begin{align}
{{u}_{r}}(r)=\frac{{{C}_{1}}Ar}{2}+\frac{{{D}_{1}}B}{\sqrt{a}}{{I}_{1}}(\sqrt{a}r)+\frac{{{C}_{2}}}{r},
\label{Solution_u}
\end{align}
{\small
in which $C_2$ is an integration constant and needs to be set to zero to keep the displacement finite at the center. Furthermore, $I_1$ is a modified Bessel function of the first kind with order 1. By substituting the boundary condition (\ref{boudcon1}) into Eq. (\ref{Solution_u}), the following relation can be obtained:
{\small
\begin{align}
\frac{{{C}_{1}}AR}{2}+\frac{{{D}_{1}}B}{\sqrt{a}}{{I}_{1}}(\sqrt{a}R)={{U}_{0}}.
\label{FBC1}
\end{align}
}
This equation contains two unknown parameters, $C_1$ and $D_1$. Consequently, an additional boundary condition must be specified from the consistent coupling boundary condition (\ref{boudcon2}). To obtain closed-form expressions for the non-zero components of the microdistortion tensor $P$, Eqs. (44) and (45) are summed together. The resultant equation is
{\small
\begin{align}
\frac{dY}{dr}+\frac{2Y}{r}=-(\frac{{{\kappa }_{\rm m}}+{{\mu }_{\rm m}}}{2{{\mu }_{\rm m}}})\frac{dZ}{dr}-\frac{{{\mu }_{\rm e}}({{\kappa }_{\rm m}}+{{\mu }_{\rm m}})}{{{\mu }_{\rm m}}({{\kappa }_{\rm e}}+{{\mu }_{\rm e}})}\frac{Z}{r}+\frac{{{C}_{1}}}{r}.
\label{dY}
\end{align}
}
Substituting Eq. (\ref{NewSol}) in Eq. (\ref{dY}) yields
{\small
\begin{align}
\frac{dY}{dr}+\frac{2Y}{r}={{D}_{1}}{{\xi }_{1}}\sqrt{a}{{I}_{1}}(\sqrt{a}r)+{{D}_{1}}{{\xi }_{2}}(\frac{{{I}_{0}}(\sqrt{a}r)}{r})+\frac{{{C}_{1}}{{\xi }_{3}}}{r},
\label{NewdY}
\end{align}
}
where
{\small
\begin{align}
{{\xi }_{1}}=-\frac{({{\kappa }_{\rm m}}+{{\mu }_{\rm m}})}{2{{\mu }_{\rm m}}}, \quad \quad {{\xi }_{2}}&=-\frac{{{\mu }_{\rm e}}({{\kappa }_{\rm m}}+{{\mu }_{\rm m}})}{{{\mu }_{\rm m}}({{\kappa }_{\rm e}}+{{\mu }_{\rm e}})}, \quad \quad {{\xi }_{3}} =\frac{{{\kappa }_{\rm m}}{{\mu }_{\rm m}}({{\kappa }_{\rm e}}+{{\mu }_{\rm e}})}{{{\kappa }_{\rm e}}{{\mu }_{\rm e}}({{\kappa }_{\rm m}}+{{\mu }_{\rm m}})+{{\kappa }_{\rm m}}{{\mu }_{\rm m}}({{\kappa }_{\rm e}}+{{\mu }_{\rm e}})}.
\label{zeta}
\end{align}
}
The solution of Eq. (\ref{NewdY}) is
{\small
\begin{align}
Y={{D}_{1}}{{\xi }_{1}}{{I}_{0}}(\sqrt{a}r)-\frac{{{D}_{1}}{{I}_{1}}(\sqrt{a}r)}{\sqrt{a}r}\left( 2{{\xi }_{1}}-{{\xi }_{2}} \right)+\frac{{{C}_{1}}{{\xi }_{3}}}{2}+\frac{{{C}_{3}}}{{{r}^{2}}},
\label{SolY}
\end{align}
{\small
where $C_3$ is an integration constant that must be set to zero in order to ensure that $Y$ remains finite at the center. Now, using boundary conditions (\ref{boudcon1}) and (\ref{boudcon2}) and Eq. (\ref{new_var}), we obtain
{\small
\begin{align}
Y(R)=0,
\label{BY}
\end{align}
}
Using Eqs. (\ref{BY}) and (\ref{SolY}), this leads to:
{\small
\begin{align}
{{D}_{1}}\left( {{\xi }_{1}}{{I}_{0}}(\sqrt{a}R)-\frac{{{I}_{1}}(\sqrt{a}R)}{\sqrt{a}R}\left( 2{{\xi }_{1}}-{{\xi }_{2}} \right) \right)+\frac{{{C}_{1}}{{\xi }_{3}}}{2}=0.
\label{FBC2}
\end{align}
}
The two unknown constants $C_1$ and $D_1$ are obtained by simultaneously solving Eqs. (\ref{FBC1}) and (\ref{FBC2}), as follows:
{\small
\begin{align}
{{C}_{1}}=\frac{2{{U}_{0}}}{R}\frac{{{I}_{1}}(\sqrt{a}R)\left( 2{{\xi }_{1}}-{{\xi }_{2}} \right)-{{\xi }_{1}}\sqrt{a}R{{I}_{0}}(\sqrt{a}R)}{A\left( {{I}_{1}}(\sqrt{a}R)\left( 2{{\xi }_{1}}-{{\xi }_{2}} \right)-{{\xi }_{1}}\sqrt{a}R{{I}_{0}}(\sqrt{a}R) \right)+B{{\xi }_{3}}{{I}_{1}}(\sqrt{a}R)},
\label{C_1}
\end{align}

\begin{align}
{{D}_{1}}=\frac{{{U}_{0}}{{\xi }_{3}}\sqrt{a}}{A\left( {{I}_{1}}(\sqrt{a}R)\left( 2{{\xi }_{1}}-{{\xi }_{2}} \right)-{{\xi }_{1}}\sqrt{a}R{{I}_{0}}(\sqrt{a}R) \right)+{{\xi }_{3}}B{{I}_{1}}(\sqrt{a}R)}.
\label{D_1}
\end{align}
}
To complete the solution procedure, it is essential to derive relations for $P_{r\theta}$ and $P_{\theta r}$. In this context, subtracting Eqs. (\ref{EAEP6}) from (\ref{EAEP5}) yields:
{\small
\begin{align}
2{{\mu }_{\text{c}}}\left( {{P}_{r\theta }}-{{P}_{\theta r}} \right)=-\frac{{{\mu }_{\text{M}}}L_{\text{c}}^{2}}{r}\left( \frac{d}{dr}\left( r\frac{d{{P}_{\theta r}}}{dr}+{{P}_{\theta r}}+{{P}_{r\theta }} \right) \right).
\label{SEq65}
\end{align}
}
By utilizing Eqs. (\ref{EAEP2}) and (\ref{EAEP5}), Eq. (\ref{SEq65}) can be rearranged as follows:
{\small
\begin{align}
\frac{d}{dr}\left( {{P}_{r\theta }}+{{P}_{\theta r}} \right)+2\left( \frac{{{P}_{r\theta }}+{{P}_{\theta r}}}{r} \right)=0.
\label{newp}
\end{align}
}
The solution of Eq. (\ref{newp}) is
{\small
\begin{align}
{{P}_{r\theta }}+{{P}_{\theta r}}=\frac{{{C}_{4}}}{{{r}^{2}}}.
\label{solnewp}
\end{align}
}
To ensure the solution remains finite at $r = 0$, $C_4$ must be set to 0 and so we have 
{\small
\begin{align}
{{P}_{r\theta }}=-{{P}_{\theta r}}.
\label{Ptr}
\end{align}
}
Substituting Eq. (\ref{Ptr}) into Eq. (\ref{EAEP2}) yields:
{\small
\begin{align}
\frac{d{{P}_{r\theta }}}{dr}=0.
\label{Prtdiff}
\end{align}
}
Equation (\ref{Prtdiff}) shows that ${{P}_{r\theta }}=C_5$  where $C_5$ is integration constant. By applying the consistent coupling boundary conditions (${{P}_{r\theta }}(R)=0$), we find ${{P}_{r\theta }}={{P}_{\theta r}}=0$. Finally, it should be noted that, in the case of the axisymmetric problem, the terms involving the Cosserat couple modulus $\mu_c$ become zero. As a result, $\mu_c$ does not appear in the solution.

Finally, the closed-form expressions for $u_r$, $P_{\theta \theta}$ and $P_{rr}$ are summarized as follows:
{\small
\begin{align}
{{u}_{r}}(r)=\frac{{{U}_{0}}}{R}\frac{A\,r\left( {{I}_{1}}(\sqrt{a}R)\left( 2{{\xi }_{1}}-{{\xi }_{2}} \right)-{{\xi }_{1}}\sqrt{a}R{{I}_{0}}(\sqrt{a}R) \right)+B{{\xi }_{3}}R{{I}_{1}}(\sqrt{a}r)}{A\left( {{I}_{1}}(\sqrt{a}R)\left( 2{{\xi }_{1}}-{{\xi }_{2}} \right)-{{\xi }_{1}}\sqrt{a}R{{I}_{0}}(\sqrt{a}R) \right)+B{{\xi }_{3}}{{I}_{1}}(\sqrt{a}R)},
\label{ur}
\end{align}

\begin{align}
{{P}_{\theta \theta }}=\frac{{{C}_{1}}}{2}\left( A-{{\xi }_{3}} \right)-{{D}_{1}}{{\xi }_{1}}{{I}_{0}}(\sqrt{a}r)+\frac{{{D}_{1}}{{I}_{1}}(\sqrt{a}r)}{\sqrt{a}r}(B+2{{\xi }_{1}}-{{\xi }_{2}}),
\label{Ptt}
\end{align}

\begin{align}
{{P}_{rr}}=Z-{{P}_{\theta \theta }}.
\label{Prr}
\end{align}
}
where $A$ is defined in Eq. (\ref{A_B}), $Z$ in Eq. (\ref{NewSol}), and ${{\xi }_{1}}$, ${{\xi }_{2}}$, and ${{\xi }_{3}}$ are given in Eq. (\ref{zeta}).

\subsection{Limit-cases} 
In this section, the limit-cases for the relaxed micromorphic model will be explored, particularly its behavior as $0 \leftarrow L_c \rightarrow \infty$. One of the simple cases occurs when both the micro and macro Poisson's ratios are zero, resulting in $\lambda_{\rm e} = \lambda_{\rm m} = 0$, which implies that $\lambda_{\rm M} = 0$. In this case, we have
{\small
\begin{align}
{{u}_{r}}(r)&={{U}_{0}}\frac{r}{R}, \\ \nonumber
{{P}_{\theta \theta }}(r)&=\frac{{{U}_{0}}}{R}\left( 1-{{\xi }_{3}} \right)+\left( \frac{\sqrt{a}r{{I}_{0}}(\sqrt{a}r)-{{I}_{1}}(\sqrt{a}r)}{\sqrt{a}R{{I}_{0}}(\sqrt{a}R)-{{I}_{1}}(\sqrt{a}R)} \right)\frac{{{U}_{0}}{{\xi }_{3}}}{r}, \\ \nonumber
P_{rr}(r)&=\frac{{{U}_{0}}}{R}\left( \frac{2{{\kappa }_{\text{e}}}{{\mu }_{\text{m}}}({{\mu }_{\text{e}}}+{{\kappa }_{\text{e}}})}{{{\mu }_{\text{e}}}{{\kappa }_{\text{e}}}({{\mu }_{\text{m}}}+{{\kappa }_{\text{m}}})+{{\mu }_{\text{m}}}{{\kappa }_{\text{m}}}({{\mu }_{\text{e}}}+{{\kappa }_{\text{e}}})}+\frac{{{\xi }_{3}}\sqrt{a}R{{I}_{0}}(\sqrt{a}r)}{\sqrt{a}R{{I}_{0}}(\sqrt{a}R)-{{I}_{1}}(\sqrt{a}R)} \right)-{{P}_{\theta \theta }},
\end{align}
}
where
{\small
\begin{align}
a=\frac{2}{L_{\text{c}}^{2}}\frac{{{\left( {{\mu }_{\text{e}}}+{{\mu }_{\text{m}}} \right)}^{2}}}{{{\mu }_{\text{e}}}{{\mu }_{\text{m}}}}, \quad \quad \quad {{\xi }_{3}}&=\frac{{{\mu }_{\text{m}}}}{{{\mu }_{\text{e}}}+{{\mu }_{\text{m}}}}.
\end{align}
}
Next, we examine two limit-cases related to classical linear elasticity. In the first scenario, it is assumed that $L_c \rightarrow 0$ represents the lower bound of macroscopic stiffness. As $L_c \rightarrow 0$, we obtain
{\small
\begin{align}
{{u}_{r}}(r)&={{U}_{0}}\frac{r}{R}, \quad \quad {{P}_{\theta \theta }}(r)=\frac{{{U}_{0}}}{R}, \quad \quad P_{rr}(r)=\frac{{{\kappa }_{\text{e}}}{{\mu }_{E}}({{\mu }_{\text{e}}}+{{\kappa }_{\text{e}}})}{{{\mu }_{\text{e}}}{{\kappa }_{\text{e}}}({{\mu }_{\text{m}}}+{{\kappa }_{\text{m}}})+{{\mu }_{\text{m}}}{{\kappa }_{\text{m}}}({{\mu }_{\text{e}}}+{{\kappa }_{\text{e}}})}\frac{2}{A}\frac{{{U}_{0}}}{R}-\frac{{{U}_{0}}}{R}.
\end{align}
}
In the second case, as $L_c$ approaches infinity, the relaxed micromorphic solution degenerates again to the known classical solution \cite{sadd2009elasticity}
{\small
\begin{align}
{{u}_{r}}(r)&={{U}_{0}}\frac{r}{R}, \quad \quad {{P}_{\theta \theta }}=\frac{{{U}_{0}}}{R}, \quad \quad P_{rr}(r)=\left( \frac{{{\kappa }_{\text{e}}}{{\mu }_{\text{m}}}({{\mu }_{\text{e}}}+{{\kappa }_{\text{e}}})}{{{\mu }_{\text{e}}}{{\kappa }_{\text{e}}}({{\mu }_{\text{m}}}+{{\kappa }_{\text{m}}})+{{\mu }_{\text{m}}}{{\kappa }_{\text{m}}}({{\mu }_{\text{e}}}+{{\kappa }_{\text{e}}})}\left( {{\xi }_{2}}-{{\xi }_{3}} \right) \right)\frac{2{{U}_{0}}}{(A{{\xi }_{2}}-B{{\xi }_{3}})R}-\frac{{{U}_{0}}}{R}.
\end{align}
}

\subsection{Numerical results and discussion}
In this section, we present numerical results to demonstrate the effects of various parameters, including the material coefficients and the characteristic length. First, variations of the non-dimensional radial displacement as a function of radial position are plotted in Figure 2 for three parameter sets as listed in Table 1, taken from Ref. \cite{sarhil2023size,sarhil2024computational}, with $R/L_{c}=2$. The figure shows that the non-dimensional radial displacement increases nearly linearly with the radial position for all parameter sets. Furthermore, it is observed that the classical linear elasticity solution coincides with the data from parameter set 1. The reason for this is that in the parameters of set 1, the $\lambda_{{\rm m}}$ and $\mu_{\rm m}$ are proportional to the macro parameters. As a result, ${\lambda_{\rm e}}{\mu_{\rm m}} = {\lambda_{\rm m}}{\mu_{\rm e}}$, leading to the outcome where $A = 1$ and $B = 0$.  In contrast, set 2 and set 3 exhibit deviations from each other and from the classical theory at various radial positions. This comparison highlights the influence of different parameter sets on the radial displacement.

\begin{table} [hp]
\caption{Values of the parameters of the relaxed micromorphic model}
\label{tab:1}
\begin{center}
\begin{tabular}{lccccc}
\hline
Set & $\lambda_{\rm M}$ (GPa) & $\mu_{\rm M}$  (GPa)& $\lambda_{\rm m}$ (GPa) & $\mu_{\rm m}$ (Gpa) & Ref. \\ \hline
1   &        17.61        &   16.13    &     30.82     &    28.23   &    \cite{sarhil2023size}  \\
2   &         1.75       &     5.90  &      11.30    &    10.19   &    \cite{sarhil2024computational}  \\
3   &          1.75       &     5.90    &    8.22   &     10.55&  \cite{sarhil2024computational} \\ \hline 
\end{tabular}
\end{center}
\end{table}

\begin{figure}[!ht]
\centering
\includegraphics[scale=0.3]{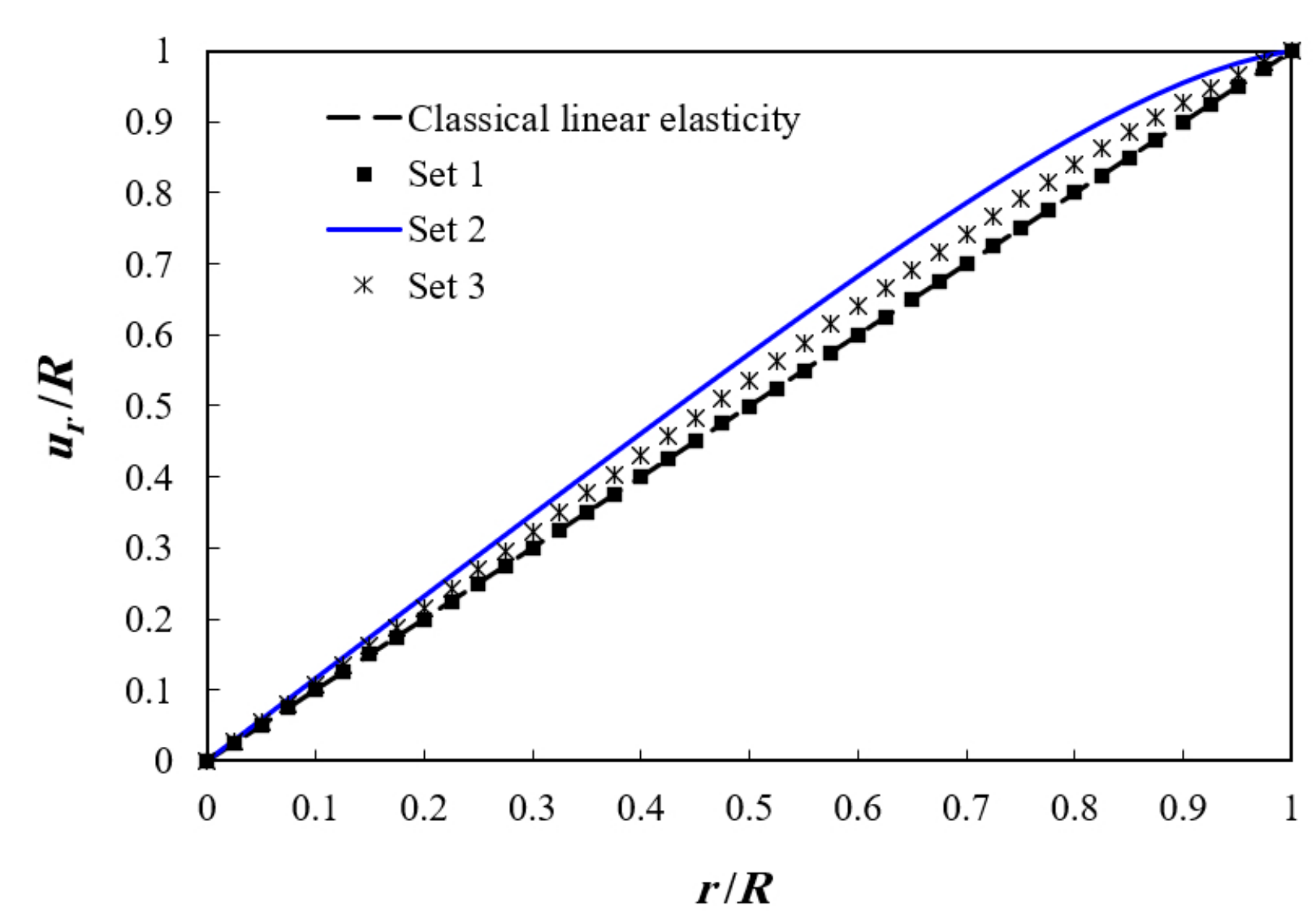}
\caption{Profile of non-dimensional radial displacement for classical elasticity and the relaxed micromorphic model for three parameter sets and $R/L_{c}=2$.}
\label{figure2}
\end{figure}

To illustrate the effect of the characteristic length, Figures 3 and 4 show the variations of the non-dimensional parameter $\delta$, defined by ${(u_{r})_{\text{Micro.}} - (u_{r})_{\text{Class.}}}/{U_0}$, as a function of radial position. These figures compare the cases where $L_{c} \rightarrow \infty$ and $L_{c} \rightarrow 0$ for the parameters of set 3. As expected, the characteristic length influences the displacement profile.
It is observed that when the characteristic length is either very large $L_{c} \rightarrow \infty$ or very small $L_{c} \rightarrow 0$, the displacement profile of the relaxed micromorphic model effectively reduces to the classical model, indicating minimal impact from the microstructure. Furthermore, the results reveal that the relaxed micromorphic model predicts a smaller displacement than the classical model at certain radial positions. This behavior highlights the critical importance of considering characteristic length in material design and analysis, as it directly influences the accuracy of displacement predictions and the overall performance of the material.

\begin{figure}[!ht]
\centering
\includegraphics[scale=0.3]{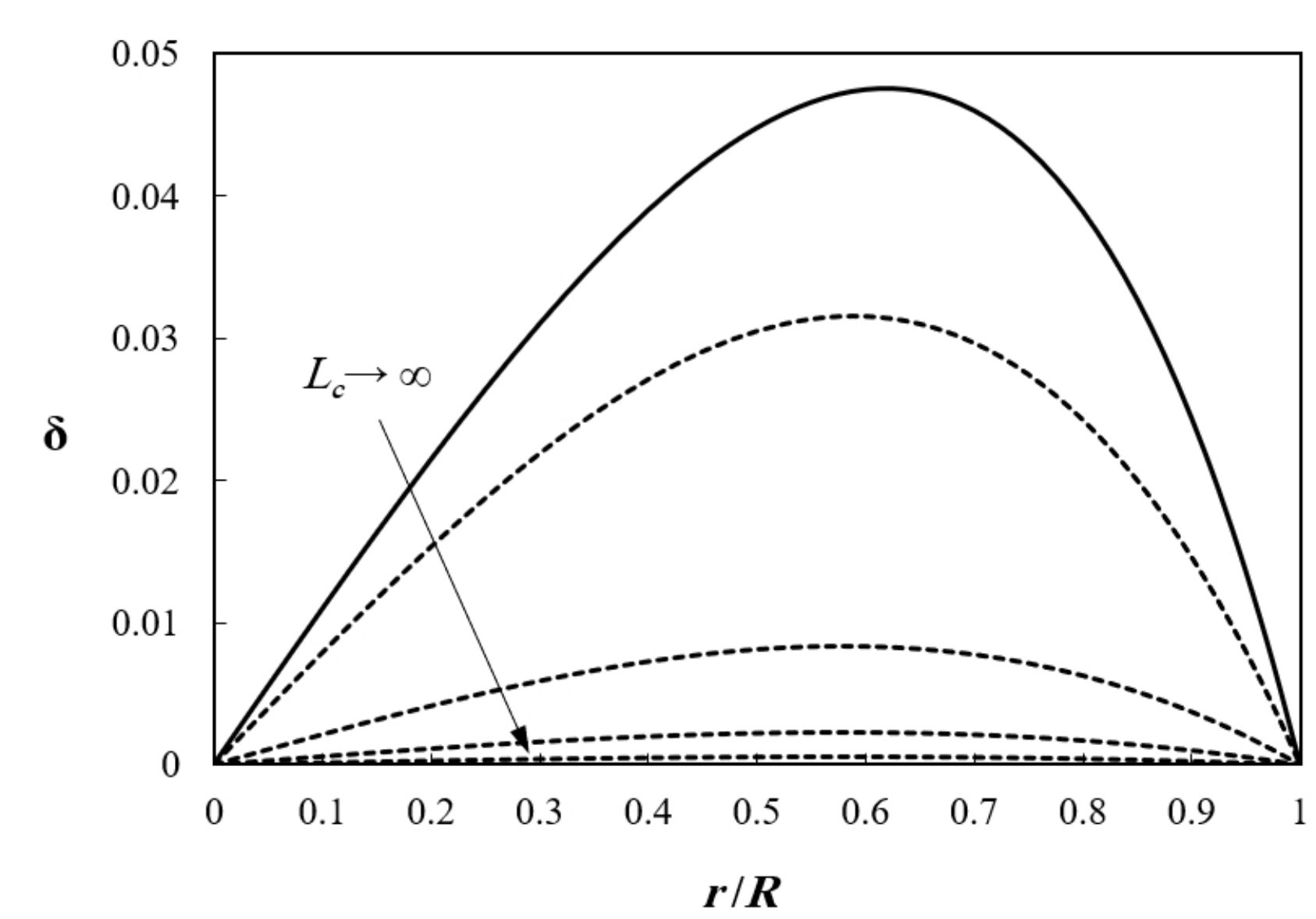}
\caption{Comparison between the relaxed micromorphic and the classical models for $R/L_{c}=\{{0.05,0.1,0.2,0.5,1}\}$.}
\label{figure3}
\end{figure}

\begin{figure}[!ht]
\centering
\includegraphics[scale=0.3]{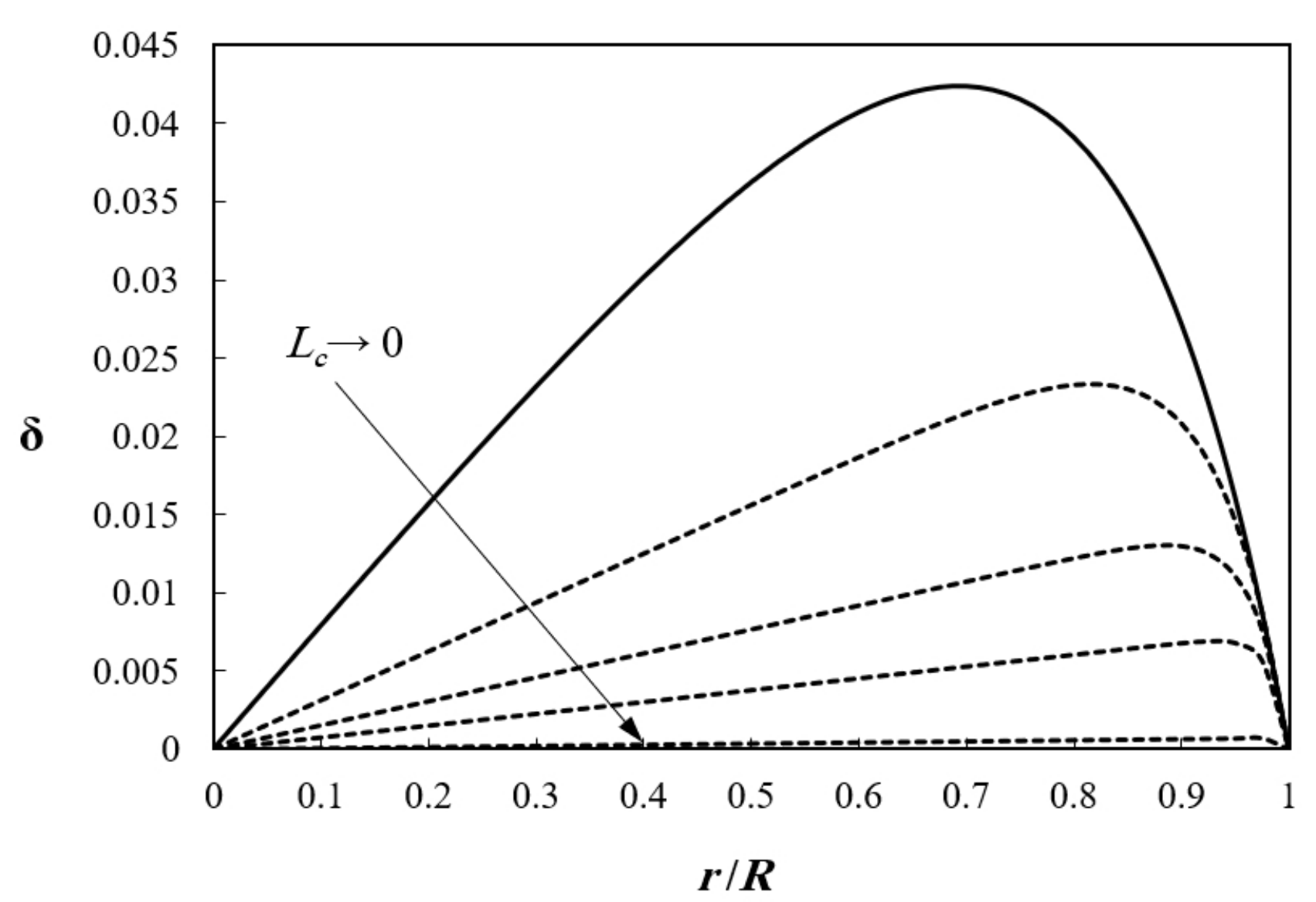}
\caption{Comparison between the relaxed micromorphic and the classical models for $R/L_{c}=\{{2,5,10,20,200}\}$.}
\label{figure4}
\end{figure}

Furthermore, to illustrate the effects of $\lambda_{\rm m}$ and $\mu_{\rm m}$, Figures 5 and 6 present the variations of the parameter $\delta$ as a function of radial position. In Figure 5, the parameters are set as $\lambda_{\rm M} = 1.75$, $\mu_{\rm M} = 5.9$, $\mu_{\rm m} = 10.55$, and $R/L_{c} = 5$, with $\lambda_{\rm m} = \beta_{1}\lambda_{\rm M}$. Similarly, in Figure 6, the parameters are $\lambda_{\rm M} = 1.75$, $\mu_{\rm M} = 5.9$, $\lambda_{\rm m} = 8.22$, and $R/L_{c} = 5$, with $\mu_{\rm m} = \beta_{2}\mu_{\rm M}$. It is observed that, for $\beta_{1} < 2$, the relaxed micromorphic model predicts a lower displacement than the classical model. Conversely, for $\beta_{1} > 2$, the behavior contrasts, indicating a complex interaction between the microstructural parameters and the overall displacement response. These findings emphasize the importance of parameter selection in the relaxed micromorphic model. Furthermore, it is observed that with increasing $\beta_{2}$ the relaxed micromorphic model reduces to the classical model.  

\begin{figure}[!ht]
\centering
\includegraphics[scale=0.31]{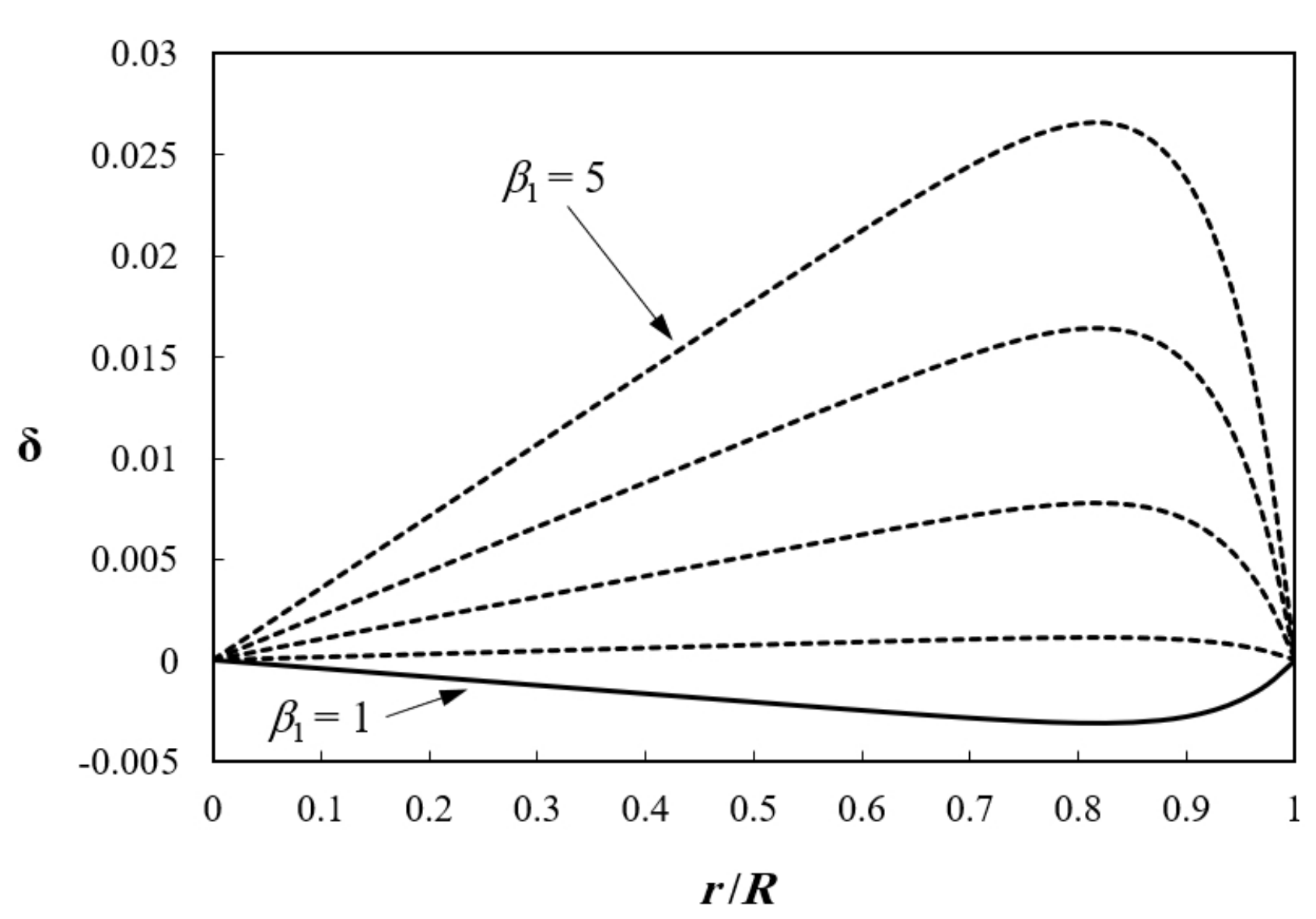}
\caption{Variations of parameter $\delta$ as a function of radial position for $\beta_{1}=\{1,2,3,4,5\}$.}
\label{figure5}
\end{figure}

\begin{figure}[!ht]
\centering
\includegraphics[scale=0.31]{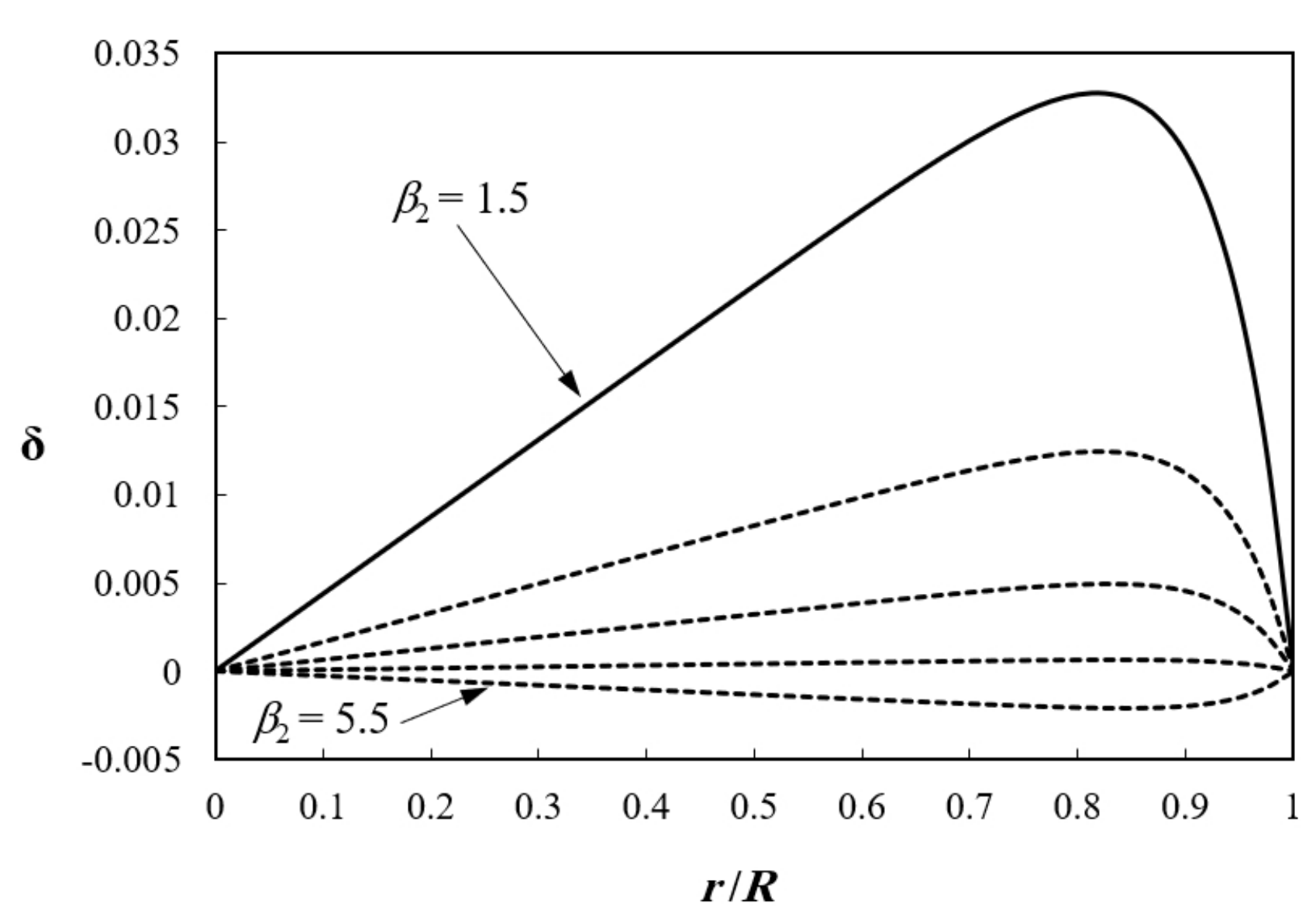}
\caption{Variations of parameter $\delta$ as a function of radial position for $\beta_{2}=\{1.5,2.5,3.5,4.5,5.5\}$.}
\label{figure6}
\end{figure}

\section{Conclusions}
In this paper, we have successfully derived the governing equations of  the isotropic relaxed micromorphic model in polar coordinates and then used the relations to solve an elastostatic axisymmetric extension problem. By utilizing modified Bessel functions, we derived closed-form solutions for both displacement and microdistortion fields using the consistent coupling boundary conditions. Our analysis demonstrated that the relaxed micromorphic model effectively captures size-dependent behaviors, which are not accounted for by classical continuum theories. Furthermore, we showed that the classical linear elasticity model can be obtained as limit-cases of the relaxed micromorphic model, highlighting its versatility and robustness. The numerical results provided insights into the effects of various parameters, such as the characteristic length and material coefficients, on the displacement. These findings underscore the importance of considering microstructural effects in material modeling to achieve accurate predictions of mechanical behavior.

\section*{Acknowledgement}
The authors are grateful to P. Gourgiotis (National Technical University of Athens), A. Sky (University of Luxembourg) and M. Sarhil (Technical University Dortmund) for critical proof reading and helpful remarks.

\begingroup
\setstretch{1}
\setlength\bibitemsep{3pt}
\printbibliography
\endgroup

\end{document}